\documentclass[a4paper,12pt]{article}
\usepackage{amssymb,latexsym,amsmath,amsthm,amscd,graphpap}
\usepackage{amsmath}
\usepackage{amsfonts}
\usepackage{amssymb}
\usepackage[dvips]{epsfig}
\newtheorem{theorem}{Theorem}[section]
\newtheorem{prop}[theorem]{Proposition}
\newtheorem{theobis}[theorem]{``Theorem''}
\newtheorem{adapted}[theorem]{Adapted Result}
\newtheorem{defi}[theorem]{Definition}

\newenvironment{demo}{ \noindent \emph{\textbf{Proof :}}}{\hfill
$\square$\\ 

\vspace{0.4cm}

}
 
\newenvironment{rem}{ \noindent \emph{\textbf{Remark :}}}{  \bigskip }

\newcommand{\Rm}{\mathbb{R}}

\newcommand{\Lm}{\mathbb{L}}

\newcommand{\Nm}{\mathbb{N}}
\newcommand{\Hm}{\mathbb{H}}

\newcommand{\Gm}{\mathfrak{G}}
\newcommand{\Gc}{\mathcal{G}}
\newcommand{\Cm}{\mathcal{C}}
\newcommand{\Nc}{\mathcal{N}}

\newcommand{\Em}{\mathcal{E}}
\newcommand{\Am}{\mathcal{A}}

\newcommand{\grad}{\overrightarrow{\nabla}}

\newcommand{\no}{n$^{\text{o}}$}

\numberwithin{equation}{section}

\begin{document}

\title{\bf Adaptation of the generic PDE's results to the notion of prevalence} 

\author{Romain JOLY\\
{\small Institut Fourier}\\ {\small UMR 5582, Universit\'e Joseph Fourier, CNRS}\\{\small 100, rue des Maths, BP74}\\{\small F-38402 St Martin d'H\`eres, FRANCE}\\ {\small \rm  Romain.Joly@ujf-grenoble.fr} }  

\date{}

\maketitle
\vspace{1cm}

%\noindent {\bf R\'esum\'e :} Un grand nombre de r\'esultats de g\'en\'ericit\'e ont \'et\'e obtenus, en particulier concernant le comportement qualitatif des solutions d'\'equations aux d\'eriv\'ees partielles. Une nouvelle notion de ``presque toujours'', la pr\'evalence, a \'et\'e r\'ecemment d\'evelopp\'ee pour les espaces vectoriels. Cette notion est int\'eressante puisque, entre autres, en dimension finie, les ensembles pr\'evalents sont exactement les ensembles de mesure de Lebesgue pleine. Le but de cet article est d'adapter les r\'esultats g\'en\'eriques en EDP \`a la notion de pr\'evalence. En particulier, nous traiterons les cas o\`u sont utilis\'es les th\'eor\`emes de Sard-Smale ou des perturbations analytiques des param\`etres.\\

%\vspace{-0.2cm}

%\noindent {\bf Mots-clefs :} pr\'evalence, g\'en\'ericit\'e, th\'eor\`eme de Smale, th\'eor\`eme de Sard-Smale.\\

%\vspace{0.1cm}

\noindent {\bf Abstract :} Many generic results have been proved, especially concerning the qualitative behaviour of solutions of partial differential equations. Recently, a new notion of ``almost always'', the prevalence, has been developped for vectorial spaces. This notion is interesting since, for example, prevalence sets are equivalent to the full Lebesgue measure sets in finite dimensional spaces. The purpose of this article is to adapt the generic PDE's results to the notion of prevalence. In particular, we consider the cases where Sard-Smale theorems or arguments of analytic perturbations of the parameters are used. \\

\vspace{-0.2cm}


\vspace{0.1cm}

\noindent {\bf AMS classification codes (2000) :} 35B30, 35P05, 37C20, 37D05, 37D15, 47F05, 58B15.\\

\vspace{1cm}

\section{Introduction}
Many important properties of partial differential equations are not always satisfied but seem to hold except for some particular cases. For example, they may hold except for a small set of coefficients of the equation. Since these properties may be very useful, one hopes to show that they ``almost always'' hold. The problem is to give a sense to this ``almost always'' notion. The study of PDE's requires infinite-dimensional spaces. Thus, the concept of set of measure zero is not very relevant as there is no natural measure on such a space. The most standard notion of ``almost always'' is the notion of generic properties, that is of properties which hold in a generic set.
\begin{defi}
Let $X$ be a Banach space and let $U\subset X$ be an open subset. A set $G\subset U$ is {\bf generic} in $U$ if it contains a countable intersection of dense open subsets of $U$. The complement of a generic set is called {\bf meager}.
\end{defi}  
\noindent This notion satisfies some natural conditions required by an ``almost always'' concept :\\
- generic sets are dense,\\
- a set which contains a generic subset is generic,\\
- a generic subset of an open generic set of $U$ is generic in $U$,\\
- a countable intersection of generic sets is generic,\\
- the notion of genericity is invariant by translation or scalar multiplication (if $U$ is a subspace).\\
To prove the genericity in the case of a PDE's property, one chooses a family of partial differential equations depending on some pertinent parameters (the domain, the potential, the diffusion coefficients... but not a universal constant or a parameter which brings the equation outside the class of models associated with the studied phenomenum). Then, one proves that the property holds for a generic set of parameters. A large number of generic results has been obtained. To give a non-exhaustive list, we can quote the generic simplicity of solutions (see \cite{Henry}, \cite{Foias-Temam} and \cite{sauttemam}), the generic hyperbolicity of equilibria of an equation (\cite{Bru-Chow}, \cite{Henry} and \cite{Smo-Was}), the generic simplicity of a spectrum (\cite{Albert}, \cite{bruno-pola}, \cite{Henry}, \cite{RJ}, \cite{OZ1}, \cite{OZ2} and \cite{Uhlenbeck}), the genericity of the Morse-Smale property (\cite{bruno-pola}, \cite{Bruno-Raugel} and \cite{RJ}) etc.\\
The main criticism against the notion of genericity is that it does not coincide with the notion of set of full Lebesgue measure in finite-dimensional spaces, which is the most relevant ``almost always'' concept in these spaces. In fact, one can even find generic subsets of $]0,1[$ which have a Lebesgue measure equal to zero. So, the question of the relevance of the notion of genericity is difficult.\\
To solve this problem, the notion of Haar-nul sets has been introduced in \cite{Christ} and the notion of prevalence in \cite{HSY}.
\begin{defi}
Let $X$ be a Banach space. A Borel set $B$ of $X$ is said {\bf Haar-nul} if there exists a finite non-negative measure $\mu\not\equiv 0$ with compact support such that 
$$\forall x\in X,~\mu(x+B)=0~.$$
More generally, any set $B\subset X$ is said Haar-nul if it is contained in a Haar-nul Borel set.\\
Let $U$ be an open subset of $X$. A set $\Em\subset U$ is said {\bf prevalent} in $U$ if $U\setminus \Em$ is a Haar-nul set of $X$.
\end{defi}
\noindent In \cite{HSY}, it is proved that the notion of prevalence satisfies the above list of properties required for an ``almost always'' concept. Moreover, in finite-dimensional spaces, prevalent sets are exactly the sets of full Lebesgue measure (see \cite{HSY}). This is an argument for using prevalence rather than genericity. On the opposite, the genericity is a notion which is stable by homeomorphism but it is not clear that the prevalence is stable, even by diffeomorphism. Another difficulty is that the notion of prevalence is not easily adaptable to the subsets of a manifold since the invariance by translation plays an important role in its definition. We refer to \cite{Ott-Yorke} for a review on prevalence.\\
What is the best notion of ``almost always''~? Which one must we use in the class of PDE's properties~? These questions are even more troublesome when one sees that prevalence and genericity may give opposite conclusions. Indeed, one can construct open dense subsets of $[0,1]$ with Lebesgue measure as small as wanted (this construction is recalled in Proposition \ref{prop-01}). By intersecting a countable number of such open dense subsets, one obtains a generic subset of $[0,1]$ with Lebesgue measure zero. Therefore, the complement of this generic subset is prevalent since it has a full Lebesgue measure. A less artificial example is given in \cite{Kahane}. It is proved in \cite{Kahane} that for all $d\geq 1$, the set of functions $f\in\Cm^0([0,1],\Rm^d)$ such that $f([0,1])$ is the adherence of its interior is a prevalent subset of $\Cm^0([0,1],\Rm^d)$. On the other hand, since any continuous function can be approximated by a piecewise-affine function, as soon as $d\geq 2$, the interior of $f([0,1])$ is generically empty. Thus, for any $d\geq 2$, the set of functions $f\in\Cm^0([0,1],\Rm^d)$ such that $f([0,1])$ is the adherence of its interior is prevalent, but negligeable from the point of view of genericity.\\
One should keep in mind that the choice of a notion of ``almost always'' for PDE's properties is not only theorical : the studied properties have concrete applications. For example, the hyperbolicity of an equilibrium and the Morse-Smale property imply respectively the local and global stability of dynamics with respect to the perturbations of the system. Thus, they are for example needed to ensure that the dynamics observed on a numerical simulation is qualitatively the same as the real dynamics of the PDE. In the same way, the simplicity of eigenvalues is required for the application of many technics.\\
Fortunately, to our knowledge, all the results of genericity in PDE's are also results of prevalence. Indeed, the purpose of this paper is to obtain what can be summarized as follows.
\begin{theobis} 
The results of the theory of partial differential equations which show the genericity of a property also show its prevalence.
\end{theobis}

Of course, we do not pretend that our paper is exhaustive as the generic results for PDE's are quite numerous. However, we have not found examples of generic properties for PDE's, the proof of which cannot be easily adapted to the notion of prevalence. This is reassuring since one does not want that the applicability of concrete methods depends on a moral and arbitrary choice between the notions of prevalence and genericity.\\
More precisely, we know two main ways to prove genericity results for PDE's. The most usual one is the use of Sard-Smale theorems or of transversality theorems (see Section \ref{sect-SS}). The second one is to obtain the density of an open property by finding, for each PDE, an explicit perturbation which brings the PDE into the desired set (see Sections \ref{sect-ana} and \ref{sect-der}). The purpose of this paper is to show that both methods can be easily adapted to the notion of prevalence.\\

\begin{rem}
a stronger notion of ``almost always'', implying both genericity and prevalence, can be found in \cite{Kolar}. However, it is too restrictive to be useful from the PDE point of view. Indeed, closed subspaces of finite codimension are not small for this notion (and so even a point is not negligible in $\Rm^n$).
\end{rem}

\noindent {\bf Acknowledgements :} I thank Alain Haraux for having introduced me to the notion of prevalence. This article is dedicated to Professor Pavol Brunovsk\'y. I am very grateful to him for having accepted to take part to my PhD committee.\\

\section{Sard-Smale theorems}\label{sect-SS}
The purpose of this section is to adapt Smale theorem to the notion of prevalence, and, as a consequence, to adapt the different versions of Sard-Smale theorem to this notion. This will show that all the genericity results obtained by using these theorems are also prevalent results. This section is also the occasion of a short review on the different versions of Sard-Smale theorem.

\subsection{The theorems of Sard and Smale}\label{sect-SaS}
Let $M$ and $M'$ be two differentiable Banach manifolds and let $f:M\longrightarrow M'$ be a differentiable map. We say that $y\in M'$ is a regular value of $f$ if, for any $x\in M$ such that $f(x)=y$, the differential $Df(x):T_xM\longrightarrow T_yM'$ is surjective. The points of $M'$ which are not regular are said critical.\\
We recall Sard theorem, the proof of which can be found in \cite{sard}. 
\begin{theorem}\label{th-sard}
Let $U$ be an open set of $\Rm^p$ and $f:U\longrightarrow \Rm^q$ be of class $\Cm^s$ with $s>max(p-q,0)$. Then, the set of critical values of $f$ in $\Rm^q$ is of Lebesgue measure zero.
\end{theorem}
Smale theorem is a generalization of Sard theorem to infinite-dimensional spaces by using the notion of Fredholm operator (see \cite{smale}). Let $X$ and $Y$ be two Banach spaces. An operator $L:X\longrightarrow Y$ is said Fredholm if its image is closed and if the dimension of its kernel and the codimension of its image are finite. The Fredholm index of $L$ is defined by Ind($L$)=dim($Ker~L$)-codim($Im~L$). Let $M$ and $M'$ be two differentiable Banach manifolds and assume that $M$ is connected. A function $f:M\longrightarrow M'$ is a Fredholm function if for all $x\in M$, its differential $Df(x):T_xM\longrightarrow T_{f(x)}M'$ is a Fredholm operator. As $M$ is connected, the index of $Df(x)$ does not depend on $x$ and is called the index of $f$.\\ 
We can easily adapt Smale theorem to the notion of prevalence by verifying that the set of regular values of a Fredholm map is not only generic but also prevalent.
\begin{theorem}\label{th-smale}
Let $M$ be a differentiable, connected and separable Banach manifold. Let $V\subset Y$ be an open set of a separable Banach space. If $f:M\longrightarrow V$ is a Fredhlom function of class $\Cm^k$ with $k>max(Ind(f),0)$, then the set of regular values of $f$ is {\bf generic and prevalent} in $V$.
\end{theorem}
\begin{demo}
The theorem of Lindel\"of shows that, from any cover by open sets of a separable metric space, we can extract a countable cover. Since the notion of generic and prevalent sets are stable by countable intersection, it is therefore sufficient to prove the result locally around a point $x_0\in M$. Up to a change of variables, we can assume that $x_0=0$, that $f(0)=0$ and that $M$ is an open set of a Banach space $X$.\\
Let $L=Df(0)$, $L$ is a Fredholm operator. We split the spaces by setting $X=X_1\oplus Ker(L)$ and $Y=Im(L)\oplus Y_2$, with $Ker(L)$ and $Y_2$ of finite dimension. As $Df(0)$ is an isomorphism from $X_1$ onto $Im(L)$, the inverse function theorem yields the existence of a neighborhood $D_1\times D_2$ of $(0,0)$ in $X_1\times Ker(L)$ such that for any $q\in D_2$, $f_{|D_1\times \{q\}}$ is a $\Cm^1$-diffeomorphism onto its image. Notice that we can choose $D_2$ compact since $Ker(L)$ is finite-dimensional. Moreover, due to the implicit function theorem, we can assume that $f$ is given by
$$f:\left(  \begin{array}{ccc} 
D_1\times D_2 & \longrightarrow & Im(L)\times Y_2 \\
(p,q) &\longmapsto & (p,g(p,q)) \end{array}\right) $$
where $g$ is a function of class $\Cm^k$ satisfying $g(0,0)=0$ and $Dg(0,0)=0$. Notice that this change of variables is possible since it concerns $X$ and not $Y$ (we have to be careful since the notion of prevalence is a priori not stable by change of variables).\\
{\bf Openess : }let $\Em$ be the set of critical values of $f$, we claim that $\Em$ is closed. Indeed, let $y_i=(p_i,r_i)$ be a sequence of critical points converging to $y=(p,r)\in Y$. Let $x_i=(p_i,q_i)$ be a point in the preimage of $y_i$ such that the differential $Df(x_i)$ is not surjective. Due to the compactness of $D_2$, we can extract a subsequence $q_i$ converging to $q\in D_2$. The continuity of $f$ and of its differential shows that $f(p,q)=y$ is a critical value. Thus, the set of regular values of $f_{|D_1\times D_2}$ is open.\\
{\bf Prevalence : }let $B_{Y_2}(0, 1)$ be the closed ball in $Y_2$ of center $0$ and radius 1. Notice that $\{0\}\times B_{Y_2}(0, 1)\subset Im(L)\times Y_2$ is a finite-dimensional disk. We define $\mu$ to be the measure on $Y$ with support equal to the disk $\{0\}\times B_{Y_2}(0,1)$ and such that $\mu$ restricted to $\{0\}\times B_{Y_2}(0, 1)$ is the Lebesgue measure on this disk.  First, $\Em$, the set of critical values of $f$, is a Borel set since it is closed. Then, let us show that $\mu(y+\Em)=0$ for any $y=(y_1,y_2)\in Y$. By definition of the support of $\mu$, it is sufficient to consider the critical values in $\{-y_1\}\times D_2$. Moreover, due to Sard theorem, the set of critical values of 
$$\tilde f:\left(  \begin{array}{ccc} 
D_2 & \longrightarrow & Y_2 \\
q &\longmapsto & g(-y_1,q) \end{array}\right) $$
is of Lebesgue mesure zero. In addition, the second component of a critical value of $f$ belonging to $\{-y_1\}\times Y_2$ must be a critical value of $\tilde f$ since $D_p f$ is surjective from $X_1$ onto $Im(L)$ and $Df=(D_p f,D\tilde f)$. Thus, $\mu(y+\Em)=0$ and the proof is complete.\\
{\bf Genericity : }as it is prevalent, the set of regular values of $f$ is dense. Thus, since it is also open, it is generic.
\end{demo}

\begin{rem}
In fact, a stronger result has been proved in \cite{Quinn-Sard}. In \cite{Quinn-Sard}, instead of the prevalence, the $s-$conullity is considered. Somehow, this notion is to the prevalence what Hausdorff dimension is to the Lebesgue measure. If $s=0$ the $s-$conullity is more or less equivalent to the prevalence, if $s>0$ the $s-$conullity is much stronger.  
\end{rem}

%\begin{rem}
%A very similar result is proved in \cite{Quinn-Sard}. We show here that the restrictions of the critical set to finite-dimensional slices are of Lebesgue measure zero. In \cite{Quinn-Sard}, the question of the Hausdorff dimension of these restrictions is considered. Notice that a result on the Hausdorff dimension is not stronger than a result on the Lebesgue measure. Indeed, we recall that there exists subsets of $[0,1]$ with positive Lebesgue measure and Hausdorff dimension less than one. A classic example is the Smith-Volterra-Cantor set which is constructed as Cantor set but by removing at each step one quarter of the set instead of one third. Its Hausdorff dimension is $\ln(2)/\ln(8/3)$ and it has a Lesbegue measure equal to $1/2$. 
%\end{rem}

Most of the generic results for PDE's are proved by using a version of Sard-Smale theorem as the ones stated in the next sections. However, in \cite{Foias-Temam}, Foias and Temam use Smale theorem only. So, we are right now able to adapt their result.
\begin{adapted}
Let $\Omega\subset \Rm^d$ ($d=2,3$) a regular bounded open domain such that $\partial \Omega$ has only a finite number of connected components $\Gamma_i$ ($i=1,\ldots,n$). Let $\nu>0$. Let $G$ be the closure of $\{f\in\mathcal{D}(\Omega)^d,~\text{div }f=0\}$ in $\Lm^2(\Omega)^d$ and let $H$ be the set of all functions $\phi\in\Hm^{3/2}(\partial\Omega)^d$ such that $\int_{\Gamma_i} \phi.\nu=0$ for $i=1,\ldots,n$.\\
Then, there exists a prevalent open set $\Gm\in G\times H$ such that, for any functions $(f,\phi)\in\Gm$, there is at most a finite number of stationary solutions of Navier-Stokes equation, that is solutions $(u,p)\in\Hm^1(\Omega)^d\times\Lm^2(\Omega)$ of 
$$\left\{\begin{array}{rcl} -\nu \Delta u+\left(u.\grad\right) u+\grad p&=&f~\text{ on }\Omega\\ \text{div }u&=&0~\text{ on }\Omega\\ u&=&\phi~\text{ on }\partial\Omega\\ \end{array}\right. $$ 
\end{adapted}

\subsection{Sard-Smale theorem}
The Sard-Smale theorem given here, and the other versions stated below, are consequences of Smale theorem. Thus, their adaptation to the notion of prevalence directly follows from Theorem \ref{th-smale}.
\begin{theorem}\label{th1}
Let $X,Y,Z$ be three Banach spaces. Let $U\subset X$, $V\subset Y$ be two open sets, $\Phi:U\times V\longrightarrow Z$ be a map of class $\Cm^k$ ($k\geq 1$) and $z$ be a point of $Z$. \\ 
We assume that :
\begin{description}
\item{i)} $\forall (x,y)\in \Phi^{-1}(z), D_x\Phi(x,y)$ is a Fredholm 
  operator of index strictly less than $k$,
\item{ii)} $\forall (x,y)\in \Phi^{-1}(z),D\Phi(x,y)$ is surjective,
\item{iii)} $X$ and $Y$ are two separable metric spaces.
\end{description}
Then $\Theta=\{y\in V/ z$ is a regular value of $\Phi(.,y)\}$ is a {\bf generic and prevalent} subset of $V$.
\end{theorem}
\begin{demo}
The outline of the proof is the following. First, we prove that $\Phi^{-1}(z)$ is a Banach manifold of $X\times Y$. Then, if $i$ is the canonical injection of this manifold in $X\times Y$ and if $\pi_Y$ is the projection onto $Y$, we show that $\pi_Y\circ i$ is a Fredholm application from $\Phi^{-1}(z)$ into $Y$ and that its index is equal to the one of $D_x\Phi$. Moreover, $z$ is a regular value of $\Phi(.,y_0)$ if and only if it is a regular value of $\pi_Y\circ i$. Thus, the conclusion directly follows from the modified Smale theorem, that is Theorem \ref{th-smale}, applied to $\pi_Y\circ i$.\\
The details of the proof are similar to the ones of the corresponding theorem of \cite{sauttemam} (see also \cite{Quinn}). \end{demo}

\begin{rem}
A particular version of Sard-Smale theorem, where $X$, $Y$ and $Z$ are finite dimensional spaces, has already been used in \cite{HSY} to prove the prevalence of some properties.
\end{rem}

\subsection{Transversal density theorems}\label{adapt-BC}
Theorem \ref{th1} can be rewritten in a more geometric frame by using the notion of transversality. The result stated here can be found in \cite{Abraham}.
\begin{defi}
Let $X$ and $Y$ be two $\Cm^1-$manifolds and $W$ be a submanifold of $Y$. We say that $f\in\Cm^1(X,Y)$ is {\bf transversal to $W$} if for every $x\in X$, either $f(x)\not\in W$ or $f(x)\in W$ and\\
i) the inverse image $(D_x f)^{-1}(T_{f(x)}W)$ splits in $T_x X$ (i.e. it is closed and admits a closed complement in $T_x X$) and\\
ii) the image $(D_x f)(T_x X)$ contains a closed complement to $T_{f(x)} W$ in $T_{f(x)}Y$.
\end{defi}
\begin{defi}
Let $\Am$, $X$ and $Y$ be $\Cm^r-$manifolds. A map $\rho:\Am\longrightarrow \Cm^r(X,Y)$ is a {\bf $\Cm^r-$representation} if the evaluation map $ev_\rho:\Am\times X\longrightarrow Y$ defined by $ev_\rho(a,x)=\rho_a(x)$ is a $\Cm^r-$map from $\Am\times X$ to $Y$.
\end{defi}

\noindent The geometric version of Sard-Smale theorem is stated as follows.
\begin{theorem}\label{th-transver}
Let $\Am$ be an open set of a Banach space. Let $X$, $Y$ be two $\Cm^r-$manifolds, $\rho:\Am\longrightarrow \Cm^r(X,Y)$ be a $\Cm^r-$representation , $W\subset Y$ be a submanifold and let $ev_\rho:\Am\times X\longrightarrow Y$ be the evaluation map.\\
Assume that\\
i) $X$ has finite dimension $n$ and $W$ has finite codimension $q$ in $Y$,\\
ii) $\Am$ and $X$ are separable,\\
iii) $r>\max(0,n-q)$,\\
iv) $ev_\rho$ is transversal to $W$.\\
Then the set $\{a\in\Am,~\rho_a$ is transversal to $W\}$ is {\bf generic and prevalent} in $\Am$.
\end{theorem}
\begin{demo}
The outline of the proof is exactly the same as the one of Theorem \ref{th1}. We refer to \cite{Abraham}. The adaptation to the notion of prevalence directly follows from the adaptation of Smale theorem.
\end{demo}

The generic version of Theorem \ref{th-transver} has been used in \cite{Bru-Chow} by Brunovsk\'y and Chow. As a consequence, their result is adapted as follows. We recall that the Whitney topology on $\Cm^k(\Rm,\Rm)$ is the topology generated by the neighborhoods given by 
\begin{equation}\label{2-topo-whitney}
\{g\in\Cm^k(\Rm,\Rm)~/~|f^{(i)}(u)-g^{(i)}(u)|\leq\delta(u),~i=0,...,k, u\in\Rm\}~,
\end{equation}
where $f$ is any function in $\Cm^k(\Rm,\Rm)$ and $\delta$ is any positive continuous function (see \cite{GG}).
\begin{adapted}
The set of functions $f\in\Cm^2(\Rm,\Rm)$ such that every solution $u\in\Hm^1(]0,1[)$ of 
\begin{equation}\label{eq-BC}
\left\{\begin{array}{ll}u_{xx}(x)+f(u(x))=0&~~x\in]0,1[\\ u(0)=u(1)=0&\end{array}\right. 
\end{equation}
is hyperbolic, is a {generic and prevalent} subset of $\Cm^k(\Rm,\Rm)$ endowed with the Whitney topology. 
\end{adapted}
We underline that the result of \cite{Bru-Chow} has been improved in different ways. In \cite{Pola}, the generic hyperbolicity of the equilibria is also proved when the Dirichlet boundary conditions are replaced by Neumann or Robin, or even more general boundary conditions. In \cite{Pereira}, the generic hyperbolicity is shown for the equilibria of the equation $u_t=(a(x)u_x)_x+f(u)$, under some particular conditions on $a(x)$. Both papers use Sard-Smale theorem or an equivalent method and are thus easily adaptable to the notion of prevalence.\\
A different approach is used in \cite{Smo-Was}, where a shorter proof of the result of \cite{Bru-Chow} is given. We discuss this approach in the last section of this paper.

\subsection{Generalizations of Sard-Smale theorem}
Theorem \ref{th1} is the most classical Sard-Smale theorem, but other versions exist. We give some examples here. Except the last one, they are all adapted to the notion of prevalence as a direct application of the modified Smale theorem, Theorem \ref{th-smale}. Theorem \ref{th-Henry} requires a little more attention.\\

\noindent{\bf $\bullet$ The non-separable case :}\\
A first generalization of Sard-Smale theorem is required when $X$ and $Y$ are not separable. Property iii) of Theorem \ref{th1} must then be replaced by the following assumption.\\
{\it iii') $\Phi$ is proper that is, if $K$ is a compact subset of $Y$, then $\{x\in X,~\Phi(x,y)=z$ with $y\in K\}$ is relatively compact.}\\
Of course, as prevalent and generic sets are stable by countable intersection, we can go further by assuming instead of iii) :\\
{\it iii'') $\Phi$ is $\sigma-$proper that is, there exists a countable family of sets which covers $\Phi^{-1}(z)$ such that, on each set, $\Phi$ is proper.}\\
A typical example where the spaces are not separable is the key theorem of \cite{bruno-pola}, Theorem 4.c.1 (a different version of this theorem has been used in \cite{Bruno-Raugel} and in \cite{RJ}). The result of Brunovsk\'y and Pol\'a\v{c}ik stated in \cite{bruno-pola} becomes :
\begin{adapted}
Let $\Omega$ be an open bounded domain of $\Rm^d$ and let $p>d$. Let $\Gm^{MS}$ be the set of non-linearities $f\in \Cm^k(\Omega\times\Rm,\Rm)$ ($k\geq 1$) such that the equilibria of the parabolic equation
$$\left\{\begin{array}{ll} u_t(x,t)=\Delta u(x,t)+f(x,u(x,t)),&(x,t)\in\Omega\times\Rm^+\\ u(x,t)=0,& (x,t)\in \partial\Omega\times\Rm^+\\ u(x,0)\in W^{1,p}(\Omega) & \end{array}\right. $$
are all hyperbolic and such that the unstable and stable manifolds of any two equilibria intersect transversally in $W^{1,p}(\Omega)$. Then, $\Gm^{MS}$ is a {generic and prevalent} subset of $\Cm^k(\Omega\times\Rm,\Rm)$ endowed with the Whitney topology.
\end{adapted}

\noindent{\bf $\bullet$ The negative index case :}\\
One notices that if $n<p$, a regular value $z$ of $f:\Rm^n\longrightarrow \Rm^p$ is in fact a point which is not in the image of $f$. Thus, if we replace assumption i) of Theorem \ref{th1} by\\
{\it i') $\forall (x,y)\in \Phi^{-1}(z), D_x\Phi(x,y)$ is a Fredholm 
  operator of negative index,}\\
then the conclusion of Theorem \ref{th1} becomes\\
{\it $\Theta=\{y\in V/ z$ is not in the image of $\Phi(.,y)\}$ is a generic and prevalent subset of $V$.}\\
As noticed in \cite{Quinn}, it is sufficient to assume that $D_x\Phi(x,y)$ is semi-Fredholm with negative index, that is that the codimension of the image can be infinite (in this case, one says that the operator is left-Fredholm).\\

\noindent{\bf $\bullet$ The non-surjective case :}\\
In general, the main difficulty for applying Sard-Smale theorem is to prove assumption ii), that is the surjectivity of $D\Phi$. In fact, in some cases, it does even not hold. Then, one needs the stronger version of Sard-Smale theorem proved by Henry (see \cite{Henry}). 
\begin{theorem}\label{th-Henry}
Let $X,Y,Z$ be three Banach spaces. Let $U\subset X$, $V\subset Y$ be two open sets, $\Phi:U\times V\longrightarrow Z$ be a map of class $\Cm^k$ ($k\geq 1$) and $z$ be a point of $Z$. \\ 
We assume that :
\begin{description}
\item{i)} $\forall (x,y)\in \Phi^{-1}(z), D_x\Phi(x,y)$ is a semi-Fredholm 
  operator (the coimage can be infinite-dimensional) and its index is strictly less than $k$,
\item{ii)} $\forall (x,y)\in \Phi^{-1}(z)$, either $D\Phi(x,y)$ is surjective, \item{~~~}~~~or dim[$Im(D\Phi(x,y))/Im(D_x\Phi(x,y))$]~$>$ dim[$Ker(D_x\Phi(x,y))$],
\item{iii)} $\Phi$ is $\sigma-$proper.
\end{description}
Then $\Theta=\{y\in V/ z$ is a regular value of $\Phi(.,y)\}$ is a {\bf generic and prevalent} subset of $V$.\\
If moreover the index of $D_x\Phi(x,y)$ is negative, then $\Theta=\{y\in V/ z$ is not in the image of $\Phi(.,y)\}$ is a {\bf generic and prevalent} subset of $V$.
\end{theorem}
\begin{demo}
As the prevalence and the genericity are stable by countable intersection, it is sufficient to prove the theorem locally around a point $(x_0,y_0)\in \Phi^{-1}(z)$ and to assume that $\Phi$ is proper. To simplify the notations, we assume without loss of generality that $(x_0,y_0)=(0,0)$ and $z=0$. If $D\Phi(0,0)$ is surjective then Theorem \ref{th-Henry} is a simple generalization of Theorem \ref{th1} similar to the ones mentioned above. If the second part of the alternative of ii) holds, then the index of $D_x\Phi(x,y)$ must be negative for all $(x,y)$ in a neighborhood of $(0,0)$. In this case, $\Phi^{-1}(z)$ is not necessarily a submanifold of $X\times Y$ and the use of Smale theorem is not possible. However, we can directly apply arguments similar to the ones used in the proof of Theorem \ref{th-smale}.\\
We set $L=D_x\Phi(0,0)$ and $N=D_y\Phi(0,0)$. We split the spaces as follows. We set $X=X_1\oplus Ker(L)$ and $n=$dim($Ker~L$). By assumption, there exists a $(n+1)-$dimensional subspace $Y_1\subset Y$ such that $NY_1$ is a $(n+1)-$dimensional subspace of $Z$ and $NY_1\cap Im(L)=\{0\}$. We set $Y=Y_1\oplus Y_2$ and $Z=Im(L)\oplus NY_1 \oplus Z_3$, where $Z_3$ is a closed subspace of $Z$, which can be reduced to $\{0\}$ or can be infinite-dimensional. We now assume that $\Phi$ is given by 
$$\Phi:\left(\begin{array}{ccc} (X_1\times Ker(L))\times(Y_1\times Y_2)&\longrightarrow& Im(L)\times NY_1 \times Z_3\\ (x_1,x_2,y_1,y_2)&\longmapsto& (Lx_1+g_1(x,y),Ny_1+g_2(x,y),g_3(x,y))  \end{array}\right) $$
with $g_i(0,0)=0$ and $D_xg_i(0,0)=0$ for $i=1,2,3$ and $D_{y_1}g_2(0,0)=0$. We want to prove that the set $\Em=\{y\in Y~/~\exists x\in X,~\Phi(x,y)=0\}$ is meager and Haar-nul in $Y$.\\
 First, $\Em$ is closed and a fortiori a Borel set. Indeed, if a sequence $(y_n)\subset\Em$ converges to $y$, we can find a sequence $(x_n)$ such that $\Phi(x_n,y_n)=0$. As $\Phi$ is proper, we may assume that $(x_n)$ converges to $x$ and, as $\Phi$ is continuous, $\Phi(x,y)=0$.\\
Then, let $B_{Y_1}(0,1)$ be the closed ball in $Y_1$ of center $0$ and radius 1. Notice that this ball is of dimension $n+1$. We define the measure $\mu$ in $Y$ as the measure with support equal to the disk $B_{Y_1}(0,1)\oplus\{0\}$ and which is the Lebesgue measure on this disk. We want to show that $\mu(\Em+y)=0$ for every $y\in Y$. This is equivalent to show that, for any fixed element $y_2\in Y_2$, the set $\{y_1\in B_{Y_1}(0,1) ~/~\exists x\in X,~\Phi(x,y_1,y_2)=0\}$ is of Lebesgue measure zero. As $D_xg_1(0,0)=0$, the implicit function theorem shows that, locally, there exists a function $\psi\in\Cm^k(Ker(L)\times Y, X_1)$ such that $Lx_1+g_1(x,y)=0$ if and only if $x_1=\psi(x_2,y)$. Thus, it is sufficient to show that, for any $y_2\in Y_2$, the set $\{y_1\in B_{Y_1}(0,1) ~/~\exists x_2\in Ker(L),~Ny_1+ g_2(\psi(x_2,y),x_2,y_1,y_2)=0\}$ is of Lebesgue measure zero. As $D_x g_2(0,0)=D_{y_1}g_2(0,0)=0$, another use of the implicit function theorem shows that there exists a function $\varphi\in\Cm^k(Ker(L), Y_1)$ such that, when $y_2$ is fixed, $Ny_1+g_2(\psi(x_2,y),x_2,y)=0$ if and only if $y_1=\varphi(x_2)$. Thus, we are reduced to show that  $\{y_1\in B_{Y_1}(0,1) ~/~y_1$ is in the image of $\varphi\}$ is of Lebesgue measure zero. As dim($Y_2$)$>$dim($Ker~L)$ by construction, this is a consequence of Sard theorem (see Theorem \ref{th-sard}) and $\Em$ is Haar-nul. Since $\Em$ is closed and Haar-nul (and thus with no interior), $\Em$ is meager. \end{demo}

In \cite{Henry}, Theorem \ref{th-Henry} is used to obtain the genericity of many properties with respect to the domain. We only give here one example. We define the space of domains of class $\Cm^2$ in $\Rm^d$ ($d\geq 1$) as the set of all domains $\Omega$ which are $\Cm^2-$diffeomorphic to a reference domain $\Omega_0$ of class $\Cm^2$ in $\Rm^d$ and we endow the space of domain with the topology inherited from the classical topology of $\Cm^2-$diffeomorphisms. 
\begin{adapted} 
Let $f\in\Cm^2(\Omega\times\Rm\times\Rm^d,\Rm)$. For a {generic and prevalent} set of domains $\Omega$ of class $\Cm^2$ in $\Rm^d$, all the solutions $u\in\Hm^1_0(\Omega)$ of 
$$\left\{\begin{array}{l} \Delta u+f(x,u,\grad u)=0~~~~\text{in }\Omega\\ u_{|\partial \Omega}=0\end{array}\right.$$
are simple, i.e. the linearization $\Delta+f'_u(x,u,\grad u)+f'_v(x,u,\grad u).\grad $ is an isomorphism from $\Hm^2(\Omega)\cap\Hm^1_0(\Omega)$ into $\Lm^2(\Omega)$.
\end{adapted}
Notice that the case where $f$ only depends on $x$ and $u$ had been studied in \cite{sauttemam} and only requires Theorem \ref{th1}. The case where $f$ also depends on $\grad u$ requires Theorem \ref{th-Henry} because the surjectivity assumption of Theorem \ref{th1} fails.\\

\section{Problems which analytically depend on parameters}\label{sect-ana}

Genericity results for PDE are not all proved with Sard-Smale theorems. Indeed, the problem can be often reduced to show that a countable family of conditions of the type $f_n(u)\neq 0$ hold on a generic set. If the function $f_n$ is continuous, the set $\{u,~f_n(u)\neq 0\}$ is open. It is therefore natural to try to prove the density of this set by showing that, at each point $u$ where $f_n$ vanishes, there exists a direction $v$ such that $Df_n(u).v\neq 0$.\\
The articles \cite{Albert}, \cite{OZ1} and \cite{OZ2} are typical examples where this method has been used to study the generic simplicity of the eigenvalues of some particular differential operators. For example, \cite{OZ2} proves the simplicity of the spectrum of Stokes operator, generically with respect to a two-dimensional domain $\Omega$. To this end, Ortega and Zuazua, prove that, if $\lambda$ is a eigenvalue of multiplicity $h$, $\lambda$ is the intersection of $h$ branches of eigenvalues $\lambda_i(\Omega)$ which locally depend analytically on the domain $\Omega$. Then, to separate the branches of eigenvalues, they look at a direction $\zeta$ such that the differentials of the eigenvalues with respect to the domain satisfy $D_\Omega\lambda_i(\Omega).\zeta\neq D_\Omega\lambda_j(\Omega).\zeta$. Notice that in some of the problems studied in the articles \cite{OZ1} and \cite{OZ2}, one cannot always ensure the existence of such a direction $\zeta$. However, the authors show that the unique continuation properties needed for this existence generically hold.\\
The purpose of this section is to adapt to the notion of prevalence the results using this method, for example \cite{Albert}, \cite{OZ1}, \cite{OZ2} or \cite{RJ}. Of course, it is surely possible to rewrite each proof in the frame of Sard-Smale theorem. However a simplest and most general way is to show that, if $f$ is analytic, the genericity of the set $\{u,~f(u)\neq 0\}$ is equivalent to its prevalence. Notice that the analyticity assumption is not so restrictive since it is satisfied in most of the eigenvalue problems.
\begin{theorem}\label{th-analytic}
Let $U$ be an open subset of a separable Banach space $X$. Let $f:U\longrightarrow \Rm$ be continuous and assume that $f$ is an analytic function in the following sense. For all $u\in U$, there exists a neighborhood $]-\varepsilon,\varepsilon[$ of $0$, such that, for any $x\in X$ with $\|x\|\leq 1$, the real map $t \longmapsto f(u+tx)$ is analytic from $]-\varepsilon,\varepsilon[$ into $\Rm$.\\
Then, the open set $\{u\in U,~f(u)\neq 0\}$ is {\bf dense in $U$ if and only if it is prevalent in $U$}.  
\end{theorem}
\begin{demo}
The ``if'' sense is clear, let us show the other part of the equivalence. Assume that $\{u\in U,~f(u)\neq 0\}$ is an open dense subset of $U$. For the same reason as in the proof of Theorem \ref{th-smale}, it is sufficient to show the result locally. Let $u_0\in U$ and assume that $f(u_0)=0$ (the case $f(u_0)\neq 0$ is trivial). To simplify the notations, assume without any restriction that $u_0=0$. Let $r>0$ be such that the ball $B(0,r)$ is included in $U$. By density and openess, there exists a ball $B(v,\rho)$ included in $B(0,r)\cap\{u\in U,~f(u)\neq 0\}$. We split the space $X$ by setting $X=Y\oplus \Rm.v$ and define a neighborhood $\Nc$ of $0$ by $\Nc=B_Y(0,\rho)\oplus ]-r,r[.v$. We set $\mu$ to be the Lebesgue measure on the segment $[-r,r].v$. For any $y\in B_Y(0,\rho)$, the segment $\{y\}\oplus ]-r,r[.v$ intersects $B(v,\rho)$ and thus $\{u\in U,~f(u)\neq 0\}$. By analyticity, there exists at most a finite number of points in $\{y\}\oplus ]-r,r[.v$ on which $f$ vanishes, and so the measure of this set is zero. This proves that $\{u\in U,~f(u)\neq 0\}$ is prevalent in $\Nc$.
\end{demo}

\begin{rem}
The definition of a notion of analyticity for functions defined on Banach spaces is an interesting problem. Different concepts have been developped. We emphasize that the assumption of Theorem \ref{th-analytic} seems to be weaker than the notions of analyticity which are the most often used. Actually, additionnal properties are often required in general, for example more regularity than the G\^ateaux-differentiability may be assumed. For a review on the notions of analyticity, we refer to \cite{Hille-P}.
\end{rem}

We obtain as a consequence the prevalent version of \cite{OZ2}. Notice that the representation used in \cite{OZ2} for the space of the domains is different from the one used by Henry in \cite{Henry}. Let $\Omega_0$ be a domain of class $\Cm^2$ of $\Rm^2$. The space of domain is identified to $W^{3,\infty}(\Omega)$ by associating with $u\in W^{3,\infty}(\Omega)$ the domain $\Omega+u=\{ x+u(x),~x\in\Omega\}$.
\begin{adapted}
Let $\Omega$ be an open and bounded $\Cm^2-$domain of $\Rm^2$. There exists a generic and prevalent set of domains $\Omega+u$ such that the eigenvalues $\lambda$ of the Stokes problem  
$$\left\{\begin{array}{ll} \Delta v+\grad p=\lambda v~&\text{ in }\Omega+u\\ \text{div~}v=0&\text{ in }\Omega+u\\ v\in\Hm^1_0(\Omega+u)^2,~p\in\Lm^2(\Omega)&   \end{array}\right. $$
are simple.
\end{adapted}

\noindent For functions $f$ defined from a Banach space to another Banach space, we can generalize Theorem \ref{th-analytic} as follows. 
\begin{prop}
Let $X$ and $Y$ be two Banach spaces and assume that $X$ is separable. Let $U$ be an open subset of $X$. Let $f:U\longrightarrow Y$ be a continuous function which is analytic in the sense that, for all $u\in U$, there exists a neighborhood $]-\varepsilon,\varepsilon[$ of $0$, such that, for any $x\in X$ with $\|x\|\leq 1$ and any continuous linear form $l:Y\longrightarrow \Rm$, the map $t\longmapsto l(f(u+tx))$ is analytic from $]-\varepsilon,\varepsilon[$ into $\Rm$.\\
Then, the open set $\{u\in U,~f(u)\neq 0\}$ is dense in $U$ if and only if it is prevalent in $U$.  
\end{prop}
\begin{demo}
We argue as in the proof of Theorem \ref{th-analytic}. Assume that $\{u\in U,~f(u)\neq 0\}$ is dense and that $u_0$ is such that $f(u_0)=0$. By density, we can find $v$ as close to $u_0$ as wanted such that $f(v)\neq 0$. Using Hahn-Banach theorem, there exists a linear form $l$ on $Y$ such that $l(f(v))=1$ and $\sup_{y\in Y\setminus\{0\}} |l(y)|/|y|=1$. Arguing with the analytic real function $l\circ f$ exactly as in the proof of Theorem \ref{th-analytic}, we prove that $\{u\in B(u_0,r),~l(f(u))\neq 0\}$ is a prevalent set of $B(u_0,r)$ for $r$ small enough. Thus, $\{u\in B(u_0,r),~f(u)\neq 0\}$ is a prevalent set of $B(u_0,r)$. Finally, by covering $U$ by a countable set of such neighborhoods, we obtain that $\{u\in U,~f(u)\neq 0\}$ is prevalent in $U$ since a countable intersection of prevalent sets is prevalent. \end{demo} 

\noindent Notice that the analyticity assumption in Theorem \ref{th-analytic} is necessary. Indeed, a classical construction provides a counter-example when this hypothesis is omitted.
\begin{prop}\label{prop-01}
For each $\alpha\in]0,1[$, there exists a function $f\in\Cm^\infty(]0,1[,\Rm)$ such that $f^{-1}(0)$ is a closed set with no interior and of Lebesgue measure bigger than $\alpha$.
\end{prop}
\begin{demo}
Let $(r_n)$ be the list of all the rationnal numbers of $]0,1[$ and let $\beta>0$ be a constant to be fixed later. For each $n\in\Nm$, we construct by classical arguments a non-negative function $\psi_n\in\Cm^\infty(]0,1[,\Rm)$ such that the support of $\psi_n$ is exactly $[r_n-\frac\beta {n^2};r_n+\frac\beta {n^2}]\cap ]0,1[$. Let $c_n$ be a positive number such that $\sup_{x\in]0,1[}\sup_{k\leq n} c_n|\psi^{(k)}_n(x)|\leq \frac 1{n^2}$. We set $\psi=\sum_{n\geq 0}c_n \psi_n$. By construction, $\psi$ belongs to $\Cm^\infty(]0,1[,\Rm)$ and $\psi(x)=0$ if and only if $x\not\in ]r_n-\frac\beta {n^2};r_n+\frac\beta {n^2}[\cap ]0,1[$ for all $n$. Thus $f^{-1}(0)$ is a closed set with no interior and its Lebesgue measure is larger than $1-\beta \frac{\pi^2}6$. \end{demo}

\section{Methods using probes}\label{sect-der}
Of course, we do not pretend that this paper covers all the methods used to prove generic results for PDE's. Even if the previous two sections seem to cover most of them, some papers use different technics. We give here one last example.\\
In \cite{Smo-Was}, Smoller and Wasserman give a short proof of the result of \cite{Bru-Chow} which has been adapted to the notion of prevalence in Section \ref{adapt-BC}. We will not enter here into the details of their proof but only look at their way for proving genericity for a set $\Gc\in\Cm^2(\Rm_+,\Rm)$. The purpose of this section is to show that their method, up to the add of a short remark, shows in fact also the prevalence of the set $\Gc$. Therefore, the adapted result of Section \ref{adapt-BC} is actually nothing new.\\
To prove the genericity of $\Gc$, the method of \cite{Smo-Was} is the following. One writes $\Gc=\cap_{k\in\Nm}\Gc^k$ where $\Gc^k$ are open sets. One wants to prove that $\Gc^k$ is dense in $\Cm^2(\Rm_+)$. Let $f_0\in\Cm^2(\Rm_+)$, using Sard's theorem, one shows that there exist $n\in\Nm$ large enough and $\varepsilon_0$ small enough such that, the set $\{\varepsilon\in]-\varepsilon_0,\varepsilon_0[,~u\mapsto f_0(u)+\varepsilon u^n\in\Gc^k\}$ is of full Lebesgue measure in $]-\varepsilon_0,\varepsilon_0[$. This obviously shows the density of $\Gc^k$ and thus the genericity of $\Gc$.\\
The method of \cite{Smo-Was} corresponds almost perfectly to the method of {\it probes} used in \cite{HSY} to show results of prevalence.
\begin{defi}
Let $U\subset X$ be an open subset of a Banach space. A finite-dimensional subspace $P\subset X$ is called a {\it probe} in $U$ for a set $T\subset X$ if there exists a Borel set $S$ containing $U\setminus T$ such that for all $x\in X$, $\mu_P(S+x)=0$, where $\mu_P$ denotes the Lebesgue measure on $P$. 
\end{defi}
By the definition of the prevalence and by the fact that a countable intersection of prevalent sets is prevalent, we immediately get the following property.
\begin{prop}
Let $X$ be a separable Banach space. Let $T$ be a subset of $X$. Assume that for each $x\in X$, there exists $r>0$ such that $B(x,r)\cap T$ admits a probe $P_x$ in $B(x,r)$. Then $T$ is prevalent in $X$.
\end{prop}
Let us now return to the proof of genericity contained in \cite{Smo-Was}. To show that the set $\Gc^k$ is not only generic but also prevalent, it is sufficient to see that the line $\Rm.u^n$ is a probe for $\Gc^k$ in a neighborhood of $f_0$ in $\Cm^2(\Rm_+)$. Almost all the work is done in \cite{Smo-Was}, the only remark to add is that $n$ and $\varepsilon_0$ can be choosen independently of $f$ in a small neighborhood of $f_0$ in $\Cm^2(\Rm_+)$ and thus the previous proposition can be applied.

\end{document}